\documentclass[11pt]{article}
\usepackage{amsmath}
\usepackage{pxfonts}
\usepackage{yfonts}
\usepackage{dsfont}
\usepackage{marvosym}
\usepackage{graphicx}
\usepackage{relsize}

\def\bel{\begin{equation}\label}
\def\eeq{\end{equation}}
\def\ds{\displaystyle}

\def\mt{\longrightarrow}
\def\v{\vskip 1em}

\def\ve{\varepsilon}
\def\R{\mathds R}

\def\C{\mathfrak{B}}

\def\A{{\bf A}}

\def\L{{\bf L}}

\def\I{{\bf I}}
\def\II{{\bf II}}
\def\III{{\bf III}}
\def\IV{{\bf IV}}
\def\M{{\bf M}}
\def\G{{\bf G}}

\def\alpha{\alphaup}
\def\beta{\betaup}
\def\gamma{\gammaup}
\def\delta{\deltaup}
\def\xi{{\xiup}}
\def\eta{{\etaup}}
\def\tau{{\tauup}}
\def\rho{{\rhoup}}
\def\phi{{\phiup}}
\def\psi{{\psiup}}
\def\lambda{{\lambdaup}}
\def\omega{\omegaup}
\def\varphi{{\varphiup}}
\def\gamma{{\gammaup}}

\newtheorem{remark}{Remark}[section]

\begin{document}
 \[\begin{array}{cc}\hbox{\LARGE{\bf Two-weight,  two-parameter fractional integration } }
 \\ \\
 \hbox{\LARGE{\bf with one side condition} } \end{array}\]
\v
\[\hbox{Lijuan Wang~~~~~Zhiming Wang~~~~~and~~~~~Zipeng Wang}\]
\begin{abstract}~~~~~
We study a family of strong fractional integral operators defined on $\R^n\times\R^m$ whose kernels have singularity on every coordinate subspace. As a result, we prove a two-weight $\L^p\mt\L^q$-norm inequality  by allowing only one of the weights to satisfy $\A_p\times\A_p$-condition.
\end{abstract}

\section{Introduction}
 \setcounter{equation}{0}
 Let $0<\alpha<n$ and $0<\beta<m$. Define
\bel{I_alpha,beta}
\I_{\alpha\beta}f(x,y)~=~\iint_{\R^n\times\R^m} f(u,v)\left({1\over |x-u|}\right)^{n-\alpha}\left({1\over |y-v|}\right)^{m-\beta}dudv.
\eeq 
 The study of certain operators  that  commute with a multi-parameter family of dilations,  dates back to the time of Jessen, Marcinkiewicz and Zygmund. Some important works have accomplished 
by Fefferman \cite{R.Fefferman}-\cite{R.Fefferman''}, Chang and Fefferman \cite{Chang-Fefferman},
Cordoba and Fefferman \cite{Cordoba-Fefferman},     Fefferman and Stein \cite{R.Fefferman-Stein},    M\"{u}ller, Ricci and Stein \cite{M.R.S},
Journ\'{e} \cite{Journe'} and Pipher \cite{Pipher}.

 In this paper, we consider 
 \bel{Two-weight Norm Ineq}
 \left\| \omega\I_{\alpha\beta}f\right\|_{\L^q(\R^n\times\R^m)}~\leq~\C_{p~q~\alpha~\beta~\omega~\sigma}~\left\| f\sigma\right\|_{\L^p(\R^n\times\R^m)},\qquad 1<p<q<\infty
 \eeq
for which $\omega^q$, $\sigma^{-{p\over p-1}}$ are non-negative, locally integrable functions. 
 
$\diamond$ {\small Throughout, $\C$ is regarded as a generic constant depending on its subindices.}

Over the past several decades, the weighted norm inequality of  one-parameter fractional integrals $f\ast |x|^{\alpha-n}, 0<\alpha<n$ has been extensively investigated. A number of classical results were established, for example by Hardy and Littlewood \cite{Hardy-Littlewood}, Sobolev \cite{Sobolev}, Stein and Weiss \cite{Stein-Weiss}, 
Hedberg \cite{Hedberg}, Fefferman and Muckenhoupt \cite{Fefferman-Muckenhoupt}, Muckenhoupt and Wheeden \cite{Muckenhoupt-Wheeden}, Perez \cite{Perez} and Sawyer and Wheeden \cite{Sawyer-Wheeden}. The multi-parameter analogue in (\ref{Two-weight Norm Ineq}) remains largely open.

Denote $Q, P$  as cubes in $\R^n$ and $\R^m$ respectively. It is well known that
(\ref{Two-weight Norm Ineq}) implies 
\bel{A}
\begin{array}{lr}\ds
\A^{\alpha\beta}_{pq}(\omega,\sigma)~=~
\\ \ds
\sup_{Q\times P\subset\R^n\times\R^m}|Q|^{{\alpha\over n}-1}|P|^{{\beta\over m}-1}
\left\{\iint_{Q\times P} \omega^q(x,y)dxdy\right\}^{1\over q}\left\{\iint_{Q\times P}\ \left({1\over \sigma}\right)^{p\over p-1}(x,y)dxdy\right\}^{p-1\over p}<\infty.
\end{array}
\eeq

The supremum $\A^{\alpha\beta}_{pq}(\omega,\sigma)$ is called  Muckenhoupt characteristic. Conversely, $\A^{\alpha\beta}_{pq}(\omega,\sigma)<\infty$ along is not a sufficient condition. In order to imply the norm inequality in (\ref{Two-weight Norm Ineq}), one must have some essential side conditions on $\omega^q$ and $\sigma^{-{p\over p-1}}$.

Denote $2Q\subset\R^n$ to be a cube having the same center of $Q$ with a side length $|2Q|^{1\over n}=2|Q|^{1\over n}$. Vice versa for $2P\subset\R^m$.
An non-negative and locally integrable  function $w$  satisfies the  {\it rectangle reverse doubling}  condition if there is an $\ve>0$ such that 
\bel{reverse doubling}
\begin{array}{lr}\ds
\int_{Q} w\left(x,y\right)dx~\leq~2^{-\ve n} ~\int_{2Q} w\left(x,y\right)dx\qquad\hbox{\small{for every  $Q\subset\R^n$ and $a.e$ $y\in\R^m$}},
\\\\ \ds
\int_{P} w\left(x,y\right)dy~\leq~2^{-\ve m} ~\int_{2P} w\left(x,y\right)dy\qquad\hbox{\small{for every  $P\subset\R^m$ and $a.e$ $x\in\R^n$}}.
\end{array}
\eeq
By assuming (\ref{reverse doubling}) on both $\omega^q$ and $\sigma^{-{p\over p-1}}$, a recent result is obtained by Sawyer and Wang \cite{Sawyer-Wang}.

\v
{\bf Theorem A:}~~Sawyer and Wang, 2020\\
{\it Let $\I_{\alpha\beta}$ defined in (\ref{I_alpha,beta}) for $0<\alpha<n$, $0<\beta<m$. Suppose that $\omega^q$ and $\sigma^{-{p\over p-1}}$ satisfy (\ref{reverse doubling}). We have
\bel{Result A}
 \left\| \omega\I_{\alpha\beta}f\right\|_{\L^q(\R^n\times\R^m)}~\leq~\C_{p~q}~\A^{\alpha\beta}_{pq}(\omega,\sigma)~\left\| f\sigma\right\|_{\L^p(\R^n\times\R^m)},\qquad 1<p<q<\infty.
 \eeq}

In contract to {\bf Theorem A}, the norm inequality in (\ref{Two-weight Norm Ineq}) can be never controlled by $\A^{\alpha\beta}_{pq}(\omega,\sigma)$ in (\ref{A}) if only one of $\omega^q$, $\sigma^{-{p\over p-1}}$ has assigned a side condition, except for the  trivial case: $\omega^q~\hbox{or}~\sigma^{-{p\over p-1}}$ satisfies the $\A_1\times\A_1$-condition. 
 This can be seen from the counter example constructed as follows.
 
Consider 
 \bel{example}
 \omega(x,y)~=~\left(1+|x|\right)^{-n}\left(1+|y|\right)^{-m},\qquad \sigma(x,y)~=~|x|^{-{n\over q}}|y|^{-{m\over q}}
 \eeq
 and
 \bel{Formula}
 {\alpha\over n}~=~{\beta\over m}~=~{1\over p}-{1\over q},\qquad 1<p<q<\infty.
 \eeq
 Given $\omega, \sigma$ in (\ref{example})-(\ref{Formula}),  we have $\sigma^{-{p\over p-1}}\in\A_{p\over p-1}\times\A_{p\over p-1}$ which is equivalent to $\sigma^p\in\A_p\times\A_p$. Moreover, 
 \bel{A_pq}
 \begin{array}{lr}\ds
 \A^{\alpha\beta}_{pq}(\omega,\sigma)~=~\sup_{Q\times P\subset\R^n\times\R^m}
\left\{{1\over |Q||P|}\iint_{Q\times P} \omega^q(x,y)dxdy\right\}^{1\over q}\left\{{1\over |Q||P|}\iint_{Q\times P}\ \left({1\over \sigma}\right)^{p\over p-1}(x,y)dxdy\right\}^{p-1\over p} 
\\\\ \ds~~~~~~~~~~~~~~~~
~\doteq~\A_{pq}(\omega,\sigma)~<~\infty.
 \end{array}
 \eeq
 On the other hand, we have
 \bel{infinity norm}
 \left\| \omega \I_{\alpha\beta}\sigma^{-1}\right\|_{\L^p(\R^n\times\R^m)\mt\L^q(\R^n\times\R^m)}~=~\infty.
 \eeq
 Regarding details can be found 
 in p.10 of the arXiv paper by Sawyer and Wang \cite{Sawyer-Wang*}.
 
 In the next section, we state our positive result by introducing a variant of $\A_{pq}(\omega,\sigma)$ in (\ref{A_pq}).

 \section{Statement of main result}
 \setcounter{equation}{0}
 For $1<p<q<\infty$, define
\bel{A_pqM}
\begin{array}{lr}\ds
\A^\M_{pq}(\omega,\sigma)~=\sup_{Q\times P\subset\R^n\times\R^m}
\left\{{1\over |Q||P|}\iint_{Q\times P} \omega^q(x,y)dxdy\right\}^{1\over q}\left\{{1\over |Q||P|}\iint_{Q\times P}\M \left({1\over \sigma}\right)^{p\over p-1}(x,y)dxdy\right\}^{p-1\over p}
\end{array}
\eeq
where $\M$ is a two-parameter strong maximal operator defined on $\R^n\times \R^m$.
\begin{remark}
Observe that if $\sigma^{-{p\over p-1}}\in\A_1\times\A_1$, then $\A_{pq}^\M(\omega,\sigma)$ is equivalent to $\A_{pq}(\omega,\sigma)$ in (\ref{A_pq}).
\end{remark}

{\bf Theorem One} ~~{\it Let $\I_{\alpha\beta}$ defined in (\ref{I_alpha,beta}) for $0<\alpha<n$, $0<\beta<m$. Suppose $\sigma^{p}\in\A_p\times\A_p$. We have
\bel{Result One}
\begin{array}{cc}\ds
\left\| \omega\I_{\alpha\beta}f\right\|_{\L^q(\R^n\times\R^m)}~\leq~\C_{p~q}~\A^\M_{pq}(\omega,\sigma)~\left\| f\sigma\right\|_{\L^p(\R^n\times\R^m)},
\\\\ \ds
{\alpha\over n}~=~{\beta\over m}~=~{1\over p}-{1\over q},\qquad 1<p<q<\infty.
\end{array}
\eeq }

In order to prove {\bf Theorem One}, we develop a two-parameter Hedberg's framework \cite{Hedberg}. Section 3 is devoted to some preliminary estimates. We give the proof of {\bf Theorem One} in section 4.

\section{Some preliminaries}
\setcounter{equation}{0}
Let $Q\times P$ shrink to $(x,y)\in\R^n\times\R^m$ in (\ref{A_pqM}). By applying Lebesgue differentiation theorem, we find
\bel{Crucial Result}
\omega(x,y)\left\{\M\left({1\over \sigma}\right)^{p\over p-1}(x,y)\right\}^{p-1\over p}~\leq~\A^\M_{pq}\left(\omega,\sigma\right).
\eeq
Moreover, it is clear that $\ \A_{pq}(\omega,\sigma)\leq \A^\M_{pq}(\omega,\sigma) $. The same estimate shows 
\bel{weights compare}
\omega(x,y)~\leq~ \A^\M_{pq}\left(\omega,\sigma\right)~\sigma(x,y).
\eeq
Let $\M_1$ and $\M_2$  be  standard maximal  operators on $\R^n$ and $\R^m$ respectively.
Define
\bel{Gf}
\G f(x,y)~=~\left\|\sigma\M_1 f(x,\cdot)\right\|_{\L^p\left(\R^{m}\right)}\left\| \sigma\M_2 f(\cdot,y)\right\|_{\L^p\left(\R^{n}\right)}
\eeq
where
\bel{L^p norms}
\begin{array}{lr}\ds
\left\| \sigma\M_1 f(x,\cdot)\right\|_{\L^p\left(\R^{m}\right)}~=~\left\{\int_{\R^{m}} \Big(\sigma\M_1 f\Big)^p(x,v) dv\right\}^{1\over p},
\\\\ \ds
\left\| \sigma\M_2 f(\cdot,y)\right\|_{\L^p\left(\R^{n}\right)}~=~\left\{\int_{\R^{n}} \Big(\sigma\M_2 f\Big)^p(u,y) du\right\}^{1\over p}.
\end{array}
\eeq
We have
\bel{Gf L^p}
\begin{array}{lr}\ds
\left\{\iint_{\R^{n}\times\R^{m}} \Big(\G f\Big)^p(x,y) dxdy\right\}^{1\over p}
\\\\ \ds ~~~~~~~
~=~\left\{\iint_{\R^{n}\times\R^{m}} \left\{\int_{\R^{m}} \Big(\sigma\M_1 f\Big)^p(x,v) dv\right\}
\left\{\int_{\R^{n}} \Big(\sigma\M_2 f\Big)^p(u,y) du\right\} dxdy\right\}^{1\over p}
\\\\ \ds ~~~~~~~
~=~\left\{\iint_{\R^{n}\times\R^{m}}\Big(\sigma\M_1 f\Big)^p(x,v) dx dv\right\}^{1\over p}
\left\{\iint_{\R^{n}\times\R^{m}} \Big(\sigma\M_2 f\Big)^p(u,y) du dy\right\}^{1\over p}
\\\\ \ds ~~~~~~~
~\leq~\C_p~\left\{\iint_{\R^{n}\times\R^{m}}\Big(f\sigma\Big)^p(x,y)dxdy\right\}^{2\over p}
\end{array}
\eeq
provided that $\sigma^p\in\A_p\times\A_p$.

\section{Proof of Theorem One}
\setcounter{equation}{0}
Given $(x,y)\in\R^n\times\R^m$, we split the two-parameter fractional integral defined in (\ref{I_alpha,beta})
$w.r.t$ the four regions:
\bel{Regions}
\begin{array}{lr}\ds
\I_{\rho\lambda}=\left\{(u,v)\colon|x-u|\leq \rho, |y-v|\leq \lambda\right\},\qquad \II_{\rho\lambda}=\left\{(u,v)\colon|x-u|>\rho, |y-v|> \lambda\right\},
\\\\ \ds
\III_{\rho\lambda}=\left\{(u,v)\colon|x-u|\le \rho, |y-v|> \lambda\right\},\qquad
\IV_{\rho\lambda}=\left\{(u,v)\colon|x-u|> \rho, |y-v|\le \lambda\right\}
\end{array}
\eeq
where $\rho=\rho(x,y)>0$ and $\lambda=\lambda(x,y)>0$
will be determined explicitly. 

Let $\chi$  be the indicator function. Define
\bel{indicators}
\begin{array}{lr}\ds
\chi_\rho(x)~=~\chi\left(|x-u|\leq\rho\right),\qquad \chi_\rho^\dagger(x)~=~\chi\left(|x-u|>\rho\right),
\\\\ \ds
\chi_\lambda(y)~=~\chi\left(|y-v|\leq\lambda\right),\qquad \chi_\lambda^\dagger(y)~=~\chi\left(|y-v|>\lambda\right).
\end{array}
\eeq

{\bf 1.} Consider
\bel{norm1}
\iint_{\I_{\rho\lambda}} f(u,v)\left({1\over |x-u|}\right)^{n-\alpha}\left({1\over |y-v|}\right)^{m-\beta} dudv.
\eeq
Note that $|x-u|^{\alpha-n}\chi_\rho(x)$ and $|y-v|^{\beta-m}\chi_\lambda(y)$ are  radially decreasing functions. They can be approximated  by $\sum_i a_i\chi_{Q_i}$ and $\sum_j b_j\chi_{P_j}$ respectively for which  $a_i>0$, $b_j>0$  and  $Q_i\times P_j\subset\R^{n}\times\R^m$ are rectangles centered on $(x,y)$ for every $i,j=1,2,\ldots$, such that
\bel{Approx uv}
\begin{array}{lr}\ds
\sum_i a_i|Q_i|~\sim~\int_{|x-u|\leq \rho} \left({1\over|x-u|}\right)^{n-\alpha} du
 ~\lesssim~ \rho^{\alpha},
\\\\ \ds
\sum_j b_j|P_j|~\sim~\int_{|y-v|\leq \lambda} \left({1\over|y-v|}\right)^{m-\beta} dv ~\lesssim~ \lambda^{\beta}.
\end{array}
\eeq
We thus have
\bel{norm1 est}
\begin{array}{lr}\ds
\iint_{\I_{\rho\lambda}} f(u,v)\left({1\over |x-u|}\right)^{n-\alpha}\left({1\over |y-v|}\right)^{m-\beta} dudv
\\\\ \ds\qquad
~\sim~\sum_{i,j} a_ib_j|Q_i| |P_j|\left\{ {1\over |Q_i||P_j|}\iint_{Q_i\times P_j}f(u,v) dudv\right\}
\\\\ \ds\qquad
~\lesssim~\left[\sum_i a_i|Q_i|\right]\left[\sum_j b_j |P_j|\right]\M f(x,y)
\\\\ \ds \qquad
~\lesssim~ \rho^{\alpha} \lambda^{\beta}\M f(x,y).
\end{array}
\eeq

{\bf 2.} Consider
\bel{norm2}
 \iint_{\II_{\rho\lambda}} f(u,v)\left({1\over |x-u|}\right)^{n-\alpha}\left({1\over |y-v|}\right)^{m-\beta} dudv.
\eeq
By using H\"{o}lder inequality, we find
\bel{norm2 est}
\begin{array}{lr}\ds
 \iint_{\II_{\rho\lambda}} f(u,v)\left({1\over |x-u|}\right)^{n-\alpha}\left({1\over |y-v|}\right)^{m-\beta} dudv
 \\\\ \ds
 ~\leq~ \left\|f\sigma\right\|_{\L^p\left(\R^n\times\R^m\right)}
 \left\{\iint_{\II_{\rho\lambda}} \left({1\over |x-u|}\right)^{(n-\alpha)\left({p\over p-1}\right)}
 \left({1\over |y-v|}\right)^{(m-\beta)\left({p\over p-1}\right)} \left({1\over \sigma}\right)^{p\over p-1}(u,v)dudv\right\}^{p-1\over p}.
 \end{array}
\eeq
A direct computation shows 
\bel{m,n p/q>0}
\begin{array}{lr}\ds
(n - \alpha){p\over p-1} - n = {n \over p-1}-{\alpha p\over p-1}~=~\left({n\over p}-\alpha\right)\left({p\over p-1}\right)~=~{n\over q}\left({p\over p-1}\right)~>~0,
\\\\ \ds
(m - \beta){p\over p-1} - m ={m \over p-1}-{\beta p\over p-1}~=~\left({m\over p}-\beta\right)\left({p\over p-1}\right)~=~{m\over q}\left({p\over p-1}\right)~>~0.
\end{array}
\eeq
Both $|x-u|^{(\alpha-n)\left({p\over p-1}\right)}\chi_{\rho}^\dagger(x)$ and $|y-v|^{(\beta-m)\left({p\over p-1}\right)}\chi_{\lambda}^\dagger(y)$ are radially decreasing.
We approximate them by  $\sum_i a_i\chi_{Q_i}$ and $\sum_j b_j\chi_{P_j}$ respectively as step {\bf 1} of which
\bel{Approx u*v*}
\begin{array}{lr}\ds
 \sum_i a_i\left|Q_i\right|~\sim~\int_{|x-u|> \rho} \left({1\over |x-u|}\right)^{(n-\alpha)\left({p\over p-1}\right)} du ~\lesssim~ \rho^{(\alpha-n)\left(\frac{p}{p-1}\right)+n},
\\\\ \ds
\sum_j b_j\left|P_j\right|~\sim~\int_{|y-v|> \lambda} \left({1\over|y-v|}\right)^{(m-\beta)\left({p\over p-1}\right)} dv ~\lesssim~ \lambda^{(\beta-m)\left(\frac{p}{p-1}\right)+m}.
\end{array}
\eeq
From (\ref{m,n p/q>0})-(\ref{Approx u*v*}), we have
\bel{norm2 Est1}
\begin{array}{lr}\ds
\iint_{\II_{\rho\lambda}}\left({1\over |x-u|}\right)^{(n-\alpha)\left({p\over p-1}\right)}
\left({1\over |y-v|}\right)^{(m-\beta)\left({p\over p-1} \right)}\left({1\over \sigma}\right)^{p\over p-1}(u,v)dudv
\\\\ \ds\qquad
~\sim~\sum_{i,j} a_i b_j\left|Q_i\right| \left|P_j\right|\left\{ {1\over |Q_i||P_j|}\iint_{Q_i\times P_j}\left({1\over \sigma}\right)^{p\over p-1}(u,v) dudv\right\}
\\\\ \ds \qquad
~\lesssim~\left[\sum_i a_i\left|Q_i\right|\right]\left[\sum_j b_j \left|P_j\right|\right]  \left\{\M\left({1\over \sigma}\right)^{p\over p-1}\right\}(x,y)
\\\\ \ds \qquad
~\leq~\C~\A_{pq}^\M(\omega,\sigma)~\rho^{{(\alpha-n)\left({p\over p-1}\right)}+n} \lambda^{{(\beta-m)\left({p\over p-1}\right)}+m}\left({1\over \omega}\right)^{p\over p-1}(x,y),\qquad \hbox{\small{by (\ref{Crucial Result})}}.
\end{array}
\eeq
Together with (\ref{norm2 est}), we find
\bel{norm2 L^p}
\begin{array}{lr}\ds
\iint_{\II_{\rho\lambda}} f(u,v)\left({1\over |x-u|}\right)^{n-\alpha}\left({1\over |y-v|}\right)^{m-\beta} dudv
\\\\ \ds
~\leq~\C~\A_{pq}^\M(\omega,\sigma)~\rho^{\alpha-n/p}\lambda^{\beta-m/p}\left\|f\sigma\right\|_{\L^p\left(\R^n\times\R^m\right)}\omega^{-1}(x,y).
\end{array}
\eeq

{\bf 3.} Consider
\bel{norm3}
\iint_{\III_{\rho\lambda}} f(u,v)\left({1\over |x-u|}\right)^{n-\alpha}\left({1\over |y-v|}\right)^{m-\beta} dudv.
\eeq
We have
\bel{norm3 est1}
\begin{array}{lr}\ds
\iint_{\III_{\rho\lambda}} f(u,v)\left({1\over |x-u|}\right)^{n-\alpha}\left({1\over |y-v|}\right)^{m-\beta} dudv
\\\\ \ds~~~~~~~
~\lesssim~\int_{|y-v|>\lambda} \left({1\over |y-v|}\right)^{m-\beta} \left\{ \int_{|x-u|\leq \rho}
f(u,v)\left({1\over |x-u|}\right)^{n-\alpha} du\right\} dv
\\\\ \ds~~~~~~~
~\lesssim~\rho^{\alpha} \int_{|y-v|>\lambda}\M_1 f(x,v) \left({1\over |y-v|}\right)^{m-\beta} dv
\end{array}
\eeq
where the second inequality in (\ref{norm3 est1}) is carried out as step {\bf 1}.
Next,  H\"{o}lder inequality implies
\bel{norm3 est2}
 \begin{array}{lr}\ds
\rho^{\alpha}\int_{|y-v|>\lambda} \M_1 f(x,v)  \left({1\over |y-v|}\right)^{m-\beta} dv
 \\\\ \ds \qquad
 ~\leq~\rho^{\alpha} \left\{\int_{\R^{m}}  \Big(\sigma\M_1 f\Big)^p(x,v) dv\right\}^{1\over p}
 \left\{\int_{|y-v|>\lambda} \left({1\over |y-v|}\right)^{(m-\beta)\left({p\over p-1}\right)}\left({1\over \sigma}\right)^{p\over p-1}(x,v)dv\right\}^{p-1\over p}
 \\\\ \ds\qquad
 ~\leq~\C~\A_{pq}^\M(\omega,\sigma)~\rho^{\alpha} \lambda^{\beta-m/p}\left\{\int_{\R^{m}}  \Big(\sigma\M_1 f\Big)^p(x,v) dv\right\}^{1\over p}\omega^{-1}(x,y).
 \end{array}
 \eeq
The second inequality in (\ref{norm3 est2}) is carried out in analogue to step {\bf 2}.

From (\ref{norm3 est1})-(\ref{norm3 est2}), we find
\bel{norm3 L^p}
\begin{array}{lr}\ds
 \iint_{\III_{\rho\lambda}} f(u,v)\left({1\over |x-u|}\right)^{n-\alpha}\left({1\over |y-v|}\right)^{m-\beta} dudv
 \\\\ \ds
 ~\leq~\C~\A_{pq}^\M(\omega,\sigma)~\rho^{\alpha} \lambda^{\beta-m/p}\left\| \sigma\M_1 f(x,\cdot)\right\|_{\L^p\left(\R^{m}\right)}\omega^{-1}(x,y).
\end{array}
\eeq

{\bf 4.} Consider
\bel{norm4}
\iint_{\IV_{\rho\lambda}} f(u,v)\left({1\over |x-u|}\right)^{n-\alpha}\left({1\over |y-v|}\right)^{m-\beta} dudv.
\eeq
We have
\bel{norm4 est1}
\begin{array}{lr}\ds
\iint_{\IV_{\rho\lambda}} f(u,v)\left({1\over |x-u|}\right)^{n-\alpha}\left({1\over |y-v|}\right)^{m-\beta} dudv
\\\\ \ds
~\lesssim~\int_{|x-u|>\rho} \left({1\over |x-u|}\right)^{n-\alpha} \left\{ \int_{|y-v|\leq \lambda}
f(u,v)\left({1\over |y-v|}\right)^{m-\beta} dv\right\} du
\\\\ \ds
~\lesssim~\lambda^{\beta} \int_{|x-u|>\rho} \M_2 f(u,y) \left({1\over |x-u|}\right)^{n-\alpha} du
\end{array}
\eeq
where the second inequality in (\ref{norm4 est1}) is carried out as step {\bf 1}.

By using H\"{o}lder inequality, we find
\bel{norm4 est2}
 \begin{array}{lr}\ds
\lambda^{\beta}\int_{|x-u|>\rho} \M_2 f(u,y)  \left({1\over |x-u|}\right)^{n-\alpha} du
 \\\\ \ds
 ~\leq~\lambda^{\beta} \left(\int_{\R^{n}}  \Big(\sigma\M_2 f\Big)^p(u,y) du\right)^{1\over p}
 \left(\int_{|x-u|>\rho} \left({1\over |x-u|}\right)^{(n-\alpha)\left({p\over p-1}\right)}\left({1\over \sigma}\right)^{p\over p-1}(u,y)du\right)^{p-1\over p}
 \\\\ \ds
 ~\leq~\C~\A_{pq}^\M(\omega,\sigma)~\rho^{\alpha-n/p} \lambda^{\beta}\left(\int_{\R^{n}}  \Big(\sigma\M_2 f\Big)^p(u,y) du\right)^{1\over p}\omega^{-1}(x,y)
 \end{array}
 \eeq
where the second inequality in (\ref{norm4 est2}) is carried out as step {\bf 2}.

From (\ref{norm4 est1})-(\ref{norm4 est2}), we have
\bel{norm4 L^p}
\begin{array}{lr}\ds
\iint_{\IV_{\rho\lambda}} f(u,v)\left({1\over |x-u|}\right)^{n-\alpha}\left({1\over |y-v|}\right)^{m-\beta} dudv
\\\\ \ds\qquad
~\leq~\C~\A_{pq}^\M(\omega,\sigma)~ \rho^{\alpha-n/p}\lambda^{\beta}\left\|\sigma\M_2 f(\cdot,y)\right\|_{\L^p\left(\R^{n}\right)}\omega^{-1}(x,y).
\end{array}
\eeq

{\bf 5.} Recall $\G f$ defined in (\ref{Gf}). Suppose
\bel{Case1}
\G f(x,y)~\leq~\omega\M f(x,y)\left\|f\sigma\right\|_{\L^p\left(\R^{n}\times \R^{m}\right)}.
\eeq
We choose $\rho$ and $\lambda$ simultaneously satisfying
\bel{EstCase1-1,2}
{\omega\M f(x,y)\over \left\| f\sigma\right\|_{\L^p\left(\R^{n}\times\R^{m}\right)}}~=~\rho^{-n/p}\lambda^{-m/p}
,\qquad
{ \left\| \sigma\M_1 f(x,\cdot)\right\|_{\L^p\left(\R^{m}\right)}
\over \left\| \sigma\M_2 f(\cdot,y)\right\|_{\L^p\left(\R^{n}\right)} }~=~{ \rho^{-n/p}\over \lambda^{-m/p}}.
\eeq
By solving the equations in (\ref{EstCase1-1,2}), we find
\bel{a Case1}
\rho^{-n/p}~=~\left\{ {\omega\M f(x,y)\over \left\| f\sigma\right\|_{\L^p\left(\R^{n}\times\R^{m}\right)}} { \left\| \sigma\M_1 f(x,\cdot)\right\|_{\L^p\left(\R^{m}\right)}
\over \left\|\sigma\M_2 f(\cdot,y)\right\|_{\L^p\left(\R^{n}\right)} }\right\}^{1\over 2}
\eeq
and
\bel{b Case1}
\lambda^{-m/p}~=~\left\{ {\omega\M f(x,y)\over \left\| f\sigma\right\|_{\L^p\left(\R^{n}\times\R^{m}\right)}} { \left\| \sigma\M_2 f(\cdot,y)\right\|_{\L^p\left(\R^{n}\right)}
\over \left\|\sigma\M_1 f(x,\cdot)\right\|_{\L^p\left(\R^{m}\right)} }\right\}^{1\over 2}.
\eeq
On the other hand, suppose
\bel{Case2}
\G f(x,y)~>~\omega\M f(x,y)\left\|f\sigma\right\|_{\L^p\left(\R^{n}\times \R^{m}\right)}.
\eeq
We choose $\rho$ and $\lambda$ simultaneously  satisfying
\bel{EstCase2-1,2}
{\G f(x,y)\over \left\| f\sigma\right\|^2_{\L^p\left(\R^{n}\times\R^{m}\right)}}~=~\rho^{-n/p}\lambda^{-m/p}
,\qquad
{ \left\| \sigma\M_1 f(x,\cdot)\right\|_{\L^p\left(\R^{m}\right)}
\over \left\| \sigma\M_2 f(\cdot,y)\right\|_{\L^p\left(\R^{n}\right)} }~=~{ \rho^{-n/p}\over \lambda^{-m/p}}.
\eeq
By solving the equations in (\ref{EstCase2-1,2}), we find
\bel{a Case2}
\rho^{-n/p}~=~\left\{ {\G f(x,y)\over \left\| f\sigma\right\|^2_{\L^p\left(\R^{n}\times\R^{m}\right)}} { \left\| \sigma\M_1 f(x,\cdot)\right\|_{\L^p\left(\R^{m}\right)}
\over \left\| \sigma\M_2 f(\cdot,y)\right\|_{\L^p\left(\R^{n}\right)} }\right\}^{1\over 2}
\eeq
and
\bel{b Case2}
\lambda^{-m/p}~=~\left\{ {\G f(x,y)\over \left\| f\sigma\right\|^2_{\L^p\left(\R^{n}\times\R^{m}\right)}} { \left\| \sigma\M_2 f(\cdot,y)\right\|_{\L^p\left(\R^{n}\right)}
\over \left\| \sigma\M_1 f(x,\cdot)\right\|_{\L^p\left(\R^{m}\right)} }\right\}^{1\over 2}.
\eeq

{\bf 6.} In the case of (\ref{Case1}),
by inserting (\ref{a Case1})-(\ref{b Case1}) into (\ref{norm1 est}), we have
\bel{norm1 Result Case1}
\begin{array}{lr}\ds
\iint_{\I_{\rho\lambda}} f(u,v)\left({1\over |x-u|}\right)^{n-\alpha}\left({1\over |y-v|}\right)^{m-\beta} dudv
\\\\ \ds
~\lesssim~\rho^{\alpha} \lambda^{\beta} \M f(x,y)
\\\\ \ds
~\lesssim~\Big(\omega\M f\Big)^{p\over q}(x,y)\left\|f\sigma\right\|_{\L^p\left(\R^n\times\R^m\right)}^{1-{p\over q}}\omega^{-1}(x,y)
\\\\ \ds
~\leq~\C~\A_{pq}^\M(\omega,\sigma)~\Big(\sigma\M f\Big)^{p\over q}(x,y)\left\|f\sigma\right\|_{\L^p\left(\R^n\times\R^m\right)}^{1-{p\over q}}\omega^{-1}(x,y),\qquad \hbox{\small{by (\ref{weights compare})}}.
\end{array}
\eeq

By inserting (\ref{a Case1})-(\ref{b Case1}) into (\ref{norm2 L^p}), we have
\bel{norm2 L^p Result Case1}
\begin{array}{lr}\ds
 \iint_{\II_{\rho\lambda}} f(u,v)\left({1\over |x-u|}\right)^{n-\alpha}\left({1\over |y-v|}\right)^{m-\beta} dudv
\\\\ \ds
~\leq~\C~\A_{pq}^\M(\omega,\sigma)~\rho^{\alpha-n/p}\lambda^{\beta-m/p}\left\|f\sigma\right\|_{\L^p\left(\R^n\times\R^m\right)} \omega^{-1}(x,y)
\\\\ \ds
~\leq~\C~\A_{pq}^\M(\omega,\sigma)~\Big(\omega\M f\Big)^{p\over q}(x,y)\left\|f\sigma\right\|_{\L^p\left(\R^n\times\R^m\right)}^{1-{p\over q}}\omega^{-1}(x,y)
\\\\ \ds
~\leq~\C~\A_{pq}^\M(\omega,\sigma)~\Big(\sigma\M f\Big)^{p\over q}(x,y)\left\|f\sigma\right\|_{\L^p\left(\R^n\times\R^m\right)}^{1-{p\over q}}\omega^{-1}(x,y),\qquad \hbox{\small{by (\ref{weights compare})}}.
\end{array}
\eeq
By inserting (\ref{a Case1})-(\ref{b Case1}) into (\ref{norm3 L^p}), we have
\bel{norm3 L^p Result Case1}
\begin{array}{lr}\ds
 \iint_{\III_{\rho\lambda}} f(u,v)\left({1\over |x-u|}\right)^{n-\alpha}\left({1\over |y-v|}\right)^{m-\beta} dudv
 \\\\ \ds
 ~\leq~\C~\A_{pq}^\M(\omega,\sigma)~\rho^{\alpha} \lambda^{\beta-m/p}\left\|\sigma\M_1 f(x,\cdot)\right\|_{\L^p\left(\R^{m}\right)}\omega^{-1}(x,y)
 \\\\ \ds
~\leq~\C~\A_{pq}^\M(\omega,\sigma)~\left\{{\omega\M f(x,y)\over \left\|f\sigma\right\|_{\L^p(\R^n\times\R^m)}}\right\}^{p\over q}
\left\{{\G f(x,y)\over \omega \M f(x,y)\left\|f\sigma\right\|_{\L^p(\R^n\times\R^m)}}\right\}^{1\over 2}
\left\|f\sigma\right\|_{\L^p\left(\R^n\times\R^m\right)}\omega^{-1}(x,y)
 \\\\ \ds
~\leq~\C~\A_{pq}^\M(\omega,\sigma)~\Big(\omega\M f\Big)^{p\over q}(x,y)\left\|f\sigma\right\|_{\L^p\left(\R^n\times\R^m\right)}^{1-{p\over q}}\omega^{-1}(x,y)
  \\\\ \ds
 ~\leq~\C~\A_{pq}^\M(\omega,\sigma)~\Big(\sigma\M f\Big)^{p\over q}(x,y)\left\|f\sigma\right\|_{\L^p\left(\R^n\times\R^m\right)}^{1-{p\over q}}\omega^{-1}(x,y), \qquad \hbox{\small{by (\ref{weights compare})}}. 
\end{array}
\eeq
By inserting (\ref{a Case1})-(\ref{b Case1}) into (\ref{norm4 L^p}), we have
\bel{norm4 L^p Result Case1}
\begin{array}{lr}\ds
 \iint_{\IV_{\rho\lambda}} f(u,v)\left({1\over |x-u|}\right)^{n-\alpha}\left({1\over |y-v|}\right)^{m-\beta} dudv
 \\\\ \ds
  ~\leq~\C~\A_{pq}^\M(\omega,\sigma)~\rho^{\alpha-n/p}\lambda^{\beta}\left\| \sigma\M_2 f(\cdot,y)\right\|_{\L^p\left(\R^{n}\right)}\omega^{-1}(x,y)
\\\\ \ds
 ~\leq~\C~\A_{pq}^\M(\omega,\sigma)~\left\{{\omega\M f(x,y)\over \left\|f\sigma\right\|_{\L^p(\R^n\times\R^m)}}\right\}^{p\over q}
\left\{{\G f(x,y)\over \omega \M f(x,y)\left\|f\sigma\right\|_{\L^p(\R^n\times\R^m)}}\right\}^{1\over 2}
\left\|f\sigma\right\|_{\L^p\left(\R^n\times\R^m\right)}\omega^{-1}(x,y)
 \\\\ \ds
 ~\leq~\C~\A_{pq}^\M(\omega,\sigma)~\Big(\omega\M f\Big)^{p\over q}(x,y)\left\|f\sigma\right\|_{\L^p\left(\R^n\times\R^m\right)}^{1-{p\over q}}\omega^{-1}(x,y)
  \\\\ \ds
  ~\leq~\C~\A_{pq}^\M(\omega,\sigma)~\Big(\sigma\M f\Big)^{p\over q}(x,y)\left\|f\sigma\right\|_{\L^p\left(\R^n\times\R^m\right)}^{1-{p\over q}}\omega^{-1}(x,y),\qquad \hbox{\small{by (\ref{weights compare})}}. 
\end{array}
\eeq

{\bf 7.} In the case of  (\ref{Case2}),
by inserting (\ref{a Case2})-(\ref{b Case2}) into (\ref{norm1 est}), we have
\bel{norm1 Result Case2}
\begin{array}{lr}\ds
\iint_{\I_{\rho\lambda}} f(u,v)\left({1\over |x-u|}\right)^{n-\alpha}\left({1\over |y-v|}\right)^{m-\beta} dudv
~\lesssim~ \rho^{\alpha} \lambda^{\beta}\M f(x,y)
\\\\ \ds~~~~~~~~~~~~~~~~~~~~~~~~~~~~~~~~~
~\lesssim~\rho^{\alpha} \lambda^{\beta}\left\{{\G f(x,y)\over \left\|f\sigma\right\|_{\L^p\left(\R^n\times\R^m\right)}}\right\}\omega^{-1}(x,y)
\\\\ \ds~~~~~~~~~~~~~~~~~~~~~~~~~~~~~~~~~
~\lesssim~\Big(\G f\Big)^{p\over q}(x,y)\left\|f\sigma\right\|_{\L^p\left(\R^n\times\R^m\right)}^{1-{2p\over q}}\omega^{-1}(x,y). 
\end{array}
\eeq
By inserting (\ref{a Case2})-(\ref{b Case2}) into (\ref{norm2 L^p}), we have
\bel{norm2 L^p Result Case2}
\begin{array}{lr}\ds
 \iint_{\II_{\rho\lambda}} f(u,v)\left({1\over |x-u|}\right)^{n-\alpha}\left({1\over |y-v|}\right)^{m-\beta} dudv
\\\\ \ds
~\leq~\C~\A_{pq}^\M(\omega,\sigma)~\rho^{\alpha-n/p}\lambda^{\beta-m/p}\left\|f\sigma\right\|_{\L^p\left(\R^n\times\R^m\right)} \omega^{-1}(x,y)
\\\\ \ds
~\leq~\C~\A_{pq}^\M(\omega,\sigma)~\Big(\G f\Big)^{p\over q}(x,y)\left\|f\sigma\right\|_{\L^p\left(\R^n\times\R^m\right)}^{1-{2p\over q}}\omega^{-1}(x,y).
\end{array}
\eeq
By inserting (\ref{a Case2})-(\ref{b Case2}) into (\ref{norm3 L^p}), we have
\bel{norm3 L^p Result Case1}
\begin{array}{lr}\ds
 \iint_{\III_{\rho\lambda}} f(u,v)\left({1\over |x-u|}\right)^{n-\alpha}\left({1\over |y-v|}\right)^{m-\beta} dudv
 \\\\ \ds
 ~\leq~\C~\A_{pq}^\M(\omega,\sigma)~\rho^{\alpha} \lambda^{\beta-m/p}\left\|\sigma\M_1 f(x,\cdot)\right\|_{\L^p\left(\R^{m}\right)}\omega^{-1}(x,y)
\\\\ \ds
~\leq~\C~\A_{pq}^\M(\omega,\sigma)~\left({\G f(x,y)\over \left\|f\sigma\right\|_{\L^p\left(\R^n\times\R^m\right)}^2}\right)^{p\over q}\left\|f\sigma\right\|_{\L^p\left(\R^n\times\R^m\right)}\omega^{-1}(x,y)
\\\\ \ds
~\leq~\C~\A_{pq}^\M(\omega,\sigma)~\Big(\G f\Big)^{p\over q}(x,y)\left\|f\sigma\right\|_{\L^p\left(\R^n\times\R^m\right)}^{1-{2p\over q}}\omega^{-1}(x,y).
\end{array}
\eeq
By inserting
(\ref{a Case2})-(\ref{b Case2}) into (\ref{norm4 L^p}), we have
\bel{norm3 L^p Result Case1}
\begin{array}{lr}\ds
 \iint_{\IV_{\rho\lambda}} f(u,v)\left({1\over |x-u|}\right)^{n-\alpha}\left({1\over |y-v|}\right)^{m-\beta} dudv
 \\\\ \ds
 ~\leq~\C~\A_{pq}^\M(\omega,\sigma)~\rho^{\alpha-n/p}\lambda^{\beta}\left\|\sigma\M_2 f(x,\cdot)\right\|_{\L^p\left(\R^{n}\right)}\omega^{-1}(x,y)
\\\\ \ds
~\leq~\C~\A_{pq}^\M(\omega,\sigma)~\left\{{\G f(x,y)\over \left\|f\sigma\right\|_{\L^p\left(\R^n\times\R^m\right)}^2}\right\}^{p\over q}\left\|f\sigma\right\|_{\L^p\left(\R^n\times\R^m\right)}\omega^{-1}(x,y)
\\\\ \ds
~\leq~\C~\A_{pq}^\M(\omega,\sigma)~\Big(\G f\Big)^{p\over q}(x,y)\left\|f\sigma\right\|_{\L^p\left(\R^n\times\R^m\right)}^{1-{2p\over q}}\omega^{-1}(x,y).
\end{array}
\eeq

{\small Department of Mathematics, Westlake University}

 {\small wanglijuan@westlake.edu.cn~~~~~wangzhiming@westlake.edu.cn~~~~~~~wangzipeng@westlake.edu.cn}

 \end{document}